\documentclass[12pt,fleqn,twoside]{article}
\usepackage[cp1251]{inputenc}
\usepackage[english,ukrainian,russian]{babel}
\usepackage{amsmath,amsfonts,amssymb}
\usepackage{latexsym,hhline,graphics,epsfig,wrapfig}

\newcounter{theorem}

\newcounter{lemma}
\newcommand{\be}{\begin{equation}}
\newcommand{\ee}{\end{equation}}

\begin{document}

\large

\noindent {\small \textbf{УДК}\ \ 517.54}

\bigskip

A. К. Бахтин $^1$, Г. П. Бахтина $^2$, В. Е. Вьюн $^1$ 
\vskip 4mm
\begin{center}\textbf{ОБ ОДНОМ НЕРАВЕНСТВЕ В ЗАДАЧАХ О НЕНАЛЕГАЮЩИХ ОБЛАСТЯХ}\end{center}
\vskip 4mm

\noindent $^1$ Институт математики НАН Украины, Киев, Украина\\
\noindent E-mail: alexander.bahtin@yandex.ru, vvikev@mail.ru\\
\noindent $^2$ Национальный технический университет "Киевский политехнический институт"\,, Киев, Украина\\
\noindent E-mail: bakhtina\_galina@mail.ru

\vskip 4mm

A.K. Bakhtin $^1$, G.P. Bakhtina $^2$, V. E. Vjun $^1$

\vskip 4mm

\begin{center}\textbf{ON AN INEQUALITY IN THE PROBLEM OF NON-OVERLAPPING DOMAINS}\end{center}
\vskip 4mm

\noindent $^1$ Institute of Mathematics of NAS, Kiev, Ukraine\\
\noindent E-mail: alexander.bahtin@yandex.ru, vvikev@mail.ru\\
\noindent $^2$ National Technical University of Ukraine “Kiev Polytechnic Institute”\,, Kiev, Ukraine\\
\noindent E-mail: bakhtina\_galina@mail.ru

\vskip 4mm

{\small \it \noindent Paper is devoted to extremal problems in geometric function theory of complex variables associated with estimates of functionals defined on the systems of non-overlapping domains. In particular, we strengthen some known result in this subject.\hfill}
\vskip 1.5mm
 {\small \it \noindent Робота присвячена дослідженню екстремальних задач геометричної теорії функцій комплексної змінної, пов'язаних з оцінками функціоналів, заданих на системах неналягаючих областей. Зокрема, основна увага приділяється посиленню деякого відомого результату у даній тематиці.\hfill}

 \normalsize \vskip 3mm

В настоящее время задачи об экстремальном разбиении занимают
значительное место в геометрической теории функций комплексного
переменного и имеют богатую историю (см., например, [1--17]). Впервые
экстремальные разбиения рассматривались при получении оценок произведения
степеней конформных радиусов неналегающих областей. Эта тематика восходит
к статье М.А. Лаврентьева 1934 года [1] и впоследствии развивалась в работах
многих авторов (см., например, [2--17]). Следует отметить, что важным элементом
исследования экстремальных задач являются глубокие результаты теории квадратичных
дифференциалов, описывающие локальную и глобальную структуру их траекторий [3].

Пусть $\mathbb{N}$, $\mathbb{R}$ -- множество натуральных и вещественных чисел
соответственно, $\mathbb{C}$ -- комплексная плоскость,
$\overline{\mathbb{C}}=\mathbb{C}\bigcup\{\infty\}$ -- ее одноточечная
компактификация, $\mathbb{R^{+}}=(0,\infty)$. Пусть $r(B,a)$ --
внутренний радиус области $B\subset\overline{\mathbb{C}}$, относительно точки
$a\in B$ (см., например, [\ref{DYBININ-94}, с.~14; \ref{BAKHTIN-08}, с.~71; \ref{DYBININ-09}, с.~30]).

Пусть $n\in \mathbb{N}$, $n\geqslant2$. Систему точек
$A_{n}:=\left\{a_{k} \in \mathbb{C}:\, k=\overline{1,n}\right\}$
назовем \textbf{\emph{$n$-лучевой}}, если $|a_{k}|\in\mathbb{R^{+}}$
при $k=\overline{1,n}$, $0=\arg a_{1}<\arg a_{2}<\ldots<\arg
a_{n}<2\pi.$

Введем обозначения $P_{k}=P_{k}(A_{n}):=\{w: \arg a_{k}<\arg w<\arg
a_{k+1}\}$, $\theta_{k}:=\arg a_{k},\,a_{n+1}:=a_{1},\,
\theta_{n+1}:=2\pi,$ $\alpha_{k}:=\displaystyle\frac{1}{\pi}\arg
\displaystyle\frac{a_{k+1}}{a_{k}},$ $\alpha_{n+1}:=\alpha_{1},$
$k=\overline{1, n},$ $\sum\limits_{k=1}^{n}\alpha_{k}=2.$

Данная работа базируется на применении кусочно-разделяющего
преобразования, развитого в [\ref{DYBININ-88}, с.~48 -- 50; \ref{DYBININ-94}, с.~27 -- 30; \ref{DYBININ-09}, с.~120].

Целью данной работы является получение точных оценок сверху для
функционала следующего вида:
\begin{equation}\label{1a}J_{n}(\gamma)=\left[r\left(B_0,0\right)r\left(B_\infty,\infty\right)\right]^{\gamma}\prod\limits_{k=1}^n r\left(B_k,a_k\right),\end{equation}
где $\gamma\in\mathbb{R^{+}}$, $A_{n}=\{a_{k}\}_{k=1}^{n}$ --
$n$-лучевая система точек, расположенная на единичной окружности, $B_0$, $B_\infty$, $\{B_{k}\}_{k=1}^{n}$ --
совокупность неналегающих областей, $a_{k}\in B_{k}$, $k=\overline{0, n}$, $\infty\in B_\infty$.

При $\gamma=\frac{1}{2}$ и $n\geq2$ оценка для функционала (\ref{1a}) для системы неналегающих
областей была найдена В.Н. Дубининым [\ref{DYBININ-88}, с.~59] методом симметризации. Г.В. Кузьмина
[\ref{KYZMINA-01}, с.~267] усилила результат работы [\ref{DYBININ-88}] и показала, что данная оценка
справедлива при $\gamma\in\left(0,\frac{n^2}{8}\right]$, $n\geq2$.
Заметим, что при $n=2$ оценка для функционала (\ref{1a}) работы [\ref{KYZMINA-01}] в точности совпадает с
оценкой работы [\ref{DYBININ-88}]. В работе [\ref{Bakhtin}] получена оценка функционала (\ref{1a}) для $\gamma\in\left(0,\frac{3}{5}\right]$.

В данной работе получено усиленную оценку функционала (\ref{1a})
для значения $n=2$.

\textbf{Теорема 1.} \emph{Пусть $0<\gamma\leq\gamma_{2}$, $\gamma_{2}=0,65$. Тогда для любой $2$-лучевой системы точек $A_2=\{a_k\}_{k=1}^2$ такой, что $|a_k|=1$, $k\in\{1,2\}$ и любого
набора взаимно непересекающихся областей $B_0$, $B_1$, $B_2$, $B_\infty$ \textup{(}$a_{0}=0\in
B_0\subset\overline{\mathbb{C}}$, $\infty\in B_\infty\subset\overline{\mathbb{C}}$, $a_1\in B_1\subset\overline{\mathbb{C}}$, $a_2\in B_2\subset\overline{\mathbb{C}}$\,\textup{)}, справедливо неравенство
$$\left[r\left(B_0,0\right)r\left(B_\infty,\infty\right)\right]^{\gamma} r\left(B_1,a_1\right)r\left(B_2,a_2\right)\leqslant$$
\begin{equation}\label{25}\leqslant\left[r\left(\Lambda_0,0\right)r\left(\Lambda_\infty,\infty\right)\right]^{\gamma} r\left(\Lambda_1, \lambda_1\right)r\left(\Lambda_2, \lambda_2\right),\end{equation} где области $\Lambda_0$, $\Lambda_\infty$, $\Lambda_1$, $\Lambda_2$ и точки $0$, $\infty$, $\lambda_1$, $\lambda_2$ ---  круговые области и соответственно полюсы квадратичного дифференциала}
\begin{equation}\label{8kl}
Q(w)dw^2=-\frac{\gamma w^{4}+(4-2\gamma)w^2+\gamma}{w^2(w^2-1)^2}\,dw^2.
\end{equation}

\textbf{Доказательство теоремы 1.} Наши исследования основаны на применении разделяющего преобразования
(см., например, [\ref{DYBININ-88}, с.~48; \ref{DYBININ-94}, с.~27~--~30; \ref{DYBININ-09}, с.~120~--~124; \ref{BAKHTIN-08}, с.~87~--~92]).

Аналогично [\ref{BAKHTIN-08}, с.~261], рассмотрим систему функций $\zeta=\pi_k(w)=-i\left(e^{-i\theta_k}w\right)^\frac{1}{\alpha_k}$,\quad
$k\in\{1,2\}$. Пусть $\Omega_k^{(1)}$, $k\in\{1,2\}$, обозначает область
плоскости $\zeta$, полученную в результате объединения связной
компоненты множества $\pi_k(B_k\bigcap\overline{P}_k)$, содержащей
точку $\pi_k(a_k)$, со своим симметричным отражением относительно
мнимой оси. В свою очередь, через $\Omega_k^{(2)}$,
$k\in\{1,2\}$, обозначаем область плоскости $\mathbb{C}_\zeta$,
полученную в результате объединения связной компоненты множества
$\pi_k(B_{k+1}\bigcap\overline{P}_k)$, содержащей точку
$\pi_k(a_{k+1})$, со своим симметричным отражением относительно мнимой
оси, $B_{n+1}:=B_1$, $\pi_n(a_{n+1}):=\pi_n(a_1)$. Кроме того,
$\Omega_k^{(0)}$ будет обозначать область плоскости
$\mathbb{C}_\zeta$ , полученную в результате объединения связной
компоненты множества $\pi_k(B_0\bigcap\overline{P}_k)$, содержащей
точку $\zeta=0$, со своим симметричным отражением относительно
мнимой оси. Семейство $\{\Omega_k^{(\infty)}\}_{k=1}^2$ является результатом разделяющего преобразования произвольной
области $B_\infty$ относительно семейств $\{P_k\}_{k=1}^2$ и $\{\pi_k\}_{k=1}^2$ в точке $\zeta=\infty$.
Обозначим $\pi_k(a_k):=\omega_k^{(1)}$,
$\pi_k(a_{k+1}):=\omega_k^{(2)}$, $k\in\{1,2\}$,
$\pi_n(a_{n+1}):=\omega_n^{(2)}$.

Из определения функций $\pi_k$ вытекает, что
$$|\pi_k(w)-\omega_k^{(1)}|\sim\frac{1}{\alpha_k}|a_k|^{\frac{1}{\alpha_k}-1}\cdot|w-a_k|,\quad
w\rightarrow a_k,\quad w\in\overline{P_k},$$
$$|\pi_k(w)-\omega_k^{(2)}|\sim\frac{1}{\alpha_k}|a_{k+1}|^{\frac{1}{\alpha_k}-1}\cdot|w-a_{k+1}|,\quad
w\rightarrow a_{k+1},\quad w\in\overline{P_k},$$
$$|\pi_k(w)|\sim|w|^\frac{1}{\alpha_k},\quad
w\rightarrow 0,\quad w\in\overline{P_k}.$$ Тогда, используя
соответствующие результаты работ [\ref{DYBININ-88}, с.~54; \ref{DYBININ-94}, с.~29], имеем неравенства
\begin{equation}\label{4a}r\left(B_k,a_k\right)\leqslant\left[\frac{r\left(\Omega_k^{(1)},\omega_k^{(1)}\right)
\cdot
r\left(\Omega_{k-1}^{(2)},\omega_{k-1}^{(2)}\right)}{\frac{1}{\alpha_k}
|a_k|^{\frac{1}{\alpha_k}-1}\cdot\frac{1}{\alpha_{k-1}}
|a_k|^{\frac{1}{\alpha_{k-1}}-1}}\right]^\frac{1}{2},\end{equation}
$$k=1,2,\quad \Omega_0^{(2)}:=\Omega_n^{(2)},\quad \omega_0^{(2)}:=\omega_n^{(2)},$$
\begin{equation}\label{5a}r\left(B_0,0\right)\leqslant\left[\prod \limits_{k=1}^2
r^{\alpha_k^2}\left(\Omega_k^{(0)},0\right)\right]^\frac{1}{2},
\end{equation}
\begin{equation}\label{6a}r\left(B_\infty,\infty\right)\leqslant\left[\prod \limits_{k=1}^2
r^{\alpha_k^2}\left(\Omega_k^{(\infty)},\infty\right)\right]^\frac{1}{2}.\end{equation}
Условия реализации знака равенства в неравенствах (\ref{4a}) -- (\ref{6a}) полностью описаны в теореме 1.9 [\ref{DYBININ-94}, с.~29]. На основании этих соотношений получаем неравенство

$$J_{2}(\gamma)\leqslant\prod\limits_{k=1}^2\left(r\left(\Omega_k^{(0)},0\right)
r\left(\Omega_k^{(\infty)},\infty\right)\right)^{\frac{\gamma\alpha_{k}^{2}}{2}}\times$$
$$\times\left(\frac{r\left(\Omega_k^{(1)},\omega_k^{(1)}\right)
\cdot r\left(\Omega_k^{(2)},\omega_k^{(2)}\right)}{\left(\frac{1}{\alpha_k}\right)^2
(|a_k||a_{k+1}|)^{\frac{1}{\alpha_k}-1}}\right)^{\frac{1}{2}}.$$
Далее, учитывая методы работ [\ref{BAKHTIN-08}, с.~262; \ref{BAKHTIN-05}, с.~300; \ref{BAKHTIN-06}, с.~871], из последнего соотношения имеем
$$J_{2}(\gamma)\leqslant4\left(\prod\limits_{k=1}^2 \alpha_{k}\right)\times$$
$$\times\prod\limits_{k=1}^2\left\{\frac{r\left(\Omega_k^{(1)},\omega_k^{(1)}\right)
\cdot
r\left(\Omega_k^{(2)},\omega_k^{(2)}\right)}{\left(|a_k|^{\frac{1}{\alpha_k}}+|a_{k+1}|^{\frac{1}{\alpha_k}}\right)^2}
\left(r\left(\Omega_k^{(0)},0\right)
r\left(\Omega_k^{(\infty)},\infty\right)\right)^{\gamma\alpha_{k}^{2}}\right\}^\frac{1}{2},$$

$|\omega_k^{(1)}|=|a_k|^\frac{1}{\alpha_k}$,
$|\omega_k^{(2)}|=|a_{k+1}|^\frac{1}{\alpha_k}$,
$|\omega_k^{(1)}-\omega_k^{(2)}|=|a_k|^\frac{1}{\alpha_k}+|a_{k+1}|^\frac{1}{\alpha_k}$.
Каждое выражение, стоящее в фигурных скобках последнего неравенства, является значением функционала
\begin{equation}\label{8afd}K_{\tau}=\left[r\left(B_0,0\right)r\left(B_\infty,\infty\right)\right]^{\tau^2}\cdot \frac{r\left(B_1,a_1\right)r\left(B_2,a_2\right)}{|a_1-a_2|^2}\end{equation}
на системе неналегающих областей $\{\Omega_k^{(0)}, \Omega_k^{(1)}, \Omega_k^{(2)},\Omega_k^{(\infty)}\}$, и соответствующей системе точек $\{0, \omega_k^{(1)}, \omega_k^{(2)}, \infty\}$ ($k\in\{1,2\}$).
Оценка функционала (\ref{8afd}), в случае фиксированных полюсов, была найдена впервые В.Н. Дубининым [\ref{DYBININ-88}, \ref{DYBININ-1988}], затем -- Г.В. Кузьминой [\ref{KYZMINA-03}], Е.Г. Емельяновым [\ref{EMELYANOV-02}], А.Л. Таргонским [\ref{TARGONSKII-05}].

На основании теоремы 4.1.1 [\ref{BAKHTIN-08}, с.~167] и инвариантности функционала (\ref{8afd}) получаем оценку
$$K_{\tau}\leqslant\Phi(\tau),\quad \tau\geq0,$$ где
$\Phi(\tau)=\tau^{2\tau^2}\cdot|1-\tau|^{-(1-\tau)^2}\cdot(1+\tau)^{-(1+\tau)^2}$. Тогда
\begin{equation}
\label{tr}J_{2}(\gamma)\leqslant4\cdot\left(\prod\limits^2_{k=1} \alpha_k\right)\left[\prod\limits^2_{k=1}
\Phi(\tau_{k})\right]^{1/2}\leqslant\end{equation}
$$\leqslant\frac{4}{\gamma}\cdot\left[\prod\limits^2_{k=1}\left(\tau_{k}^{2\tau_{k}^2+2}\cdot|1-\tau_{k}|^{-(1-\tau_{k})^2}\cdot(1+\tau_{k})^{-(1+\tau_{k})^2}\right)\right]^{\frac{1}{2}},$$
где $\tau_{k}=\sqrt{\gamma}\cdot\alpha_{k}$, $k=\overline{1,2}$.

Рассмотрим подробнее функцию $$\Psi(x)=x^{2x^2+2}\cdot|1-x|^{-(1-x)^2}\cdot(1+x)^{-(1+x)^2}.$$

$\Psi(x)$ -- логарифмически выпуклая на промежутке $[0, x_{0}]$, где\quad $x_{0}\approx0,88441$, $\Psi(x_{0})=0,07002$. На промежутке $[0, x_{1}]$ ($x_{1}\approx0,58142$ -- точка максимума функции $\Psi(x)$, $\Psi(x_{1})\approx 0,08674$) функция возрастает от значения $\Psi(0)=0$ до $\Psi(x_{1})$, и убывает на промежутке $(x_{1}, \infty]$.

Далее, применяя к функции $\Psi(x)$ идеи работ [\ref{kovalev},\ref{Bakhtin}] и некоторые дополнительные рассуждения, мы получаем утверждение теоремы 1.

\vskip 3.5mm

\footnotesize
\begin{enumerate}
\Rus

\item\label{LAVRENTEV-34}{\it{Лаврентьев М.А.}\/} К теории конформных
отображений//Тр. физ.-мат. ин-та АН СССР. --- 1934. --- \textbf{5}. --- С. 159---245.

\item\label{GOLUZIN-66}{\it{Голузин Г.М.}\/} Геометрическая теория функций комплексного
переменного. --- Москва: Наука, 1966. --- 628 с.

\item\label{DJENKINS-62}{\it{Дженкинс Дж.А.}\/} Однолистные функции и конформные отображения. --- Москва: Изд-во иностр. лит., 1962. --- 256 с.

\item\label{DYBININ-88}{\it{Дубинин В. Н.}\/} Разделяющее преобразование областей и
задачи об экстремальном разбиении// Зап. науч. сем. Ленингр.
отд-ния Мат. ин-та АН СССР. --- 1988. --- \textbf{168}. --- С. 48---66.

\item\label{DYBININ-94}{\it{Дубинин В.Н.}\/} Метод симметризации в геометрической
теории функций комплексного переменного// Успехи мат. наук. ---
1994. --- \textbf{49}, \mbox{№ 1(295).} --- С. 3---76.

\item\label{KYZMINA-01}{\it{Кузьмина Г.В.}\/} Задачи об экстремальном разбиении римановой сферы//
Зап. науч. сем. ПОМИ. --- 2001. --- \textbf{276}. --- С. 253---275.

\item\label{BAKHTIN-08}{\it{Бахтин А. К., Бахтина Г. П., Зелинский Ю. Б.}\/} Тополого-алгебраические структуры  и геометрические методы в комплексном анализе// Праці ін-ту мат-ки НАН України. --- 2008. --- \textbf{73}. --- 308 с.

\item\label{DYBININ-09}{\it{Дубинин В.Н.}\/} Емкости конденсаторов и симметризация в
геометрической теории функций комплексного переменного// Владивосток: "Дальнаука" ДВО РАН, 2009. --- 390с.

\item\label{BAKHTIN-05}{\it{Бахтин А.К., Таргонский А.Л.}\/} Экстремальные задачи и
квадратичные дифференциалы//Нелінійні коливання. --- 2005. --- \textbf{8}, \mbox{№ 3.} --- С. 298---303.

\item\label{TARGONSKII-08}{\it{Бахтин А.К., Таргонский А.Л.}\/} О произведении внутренних радиусов неналегающих областей и открытых множеств//Доп. НАН України. -- 2008. --- \mbox{№ 5}. --- С. 7---12.

\item\label{BAKHTIN-06}{\it{Бахтин А.К.}\/} Экстремальные задачи о неналегающих областях со свободными полюсами на окружности//Укр. мат. журн. --- 2006. ---  \textbf{58}, \mbox{№ 7.} --- С. 867---886.

\item\label{KYZMINA-03}{\it{Кузьмина Г.В.}\/} Метод экстремальной метрики в задачах о максимуме
произведения степеней конформных радиусов неналегающих областей при
наличии свободных параметров// Зап. науч. сем. ПОМИ. --- 2003. ---  \textbf{302}. --- С. 52---67.

\item\label{EMELYANOV-02}{\it{Емельянов Е.Г.}\/} К задаче о максимуме произведения степеней
конформных радиусов неналегающих областей// Зап. науч. сем. ПОМИ.
--- 2002. --- \textbf{286}. --- С. 103---114.

\item\label{DYBININ-1988}{\it{Дубинин В.Н.}\/} Метод симметризации в геометрической теории функций:
Дис. … д-ра физ.-мат. наук. --- Владивосток, 1988. -- 193 с.

\item\label{TARGONSKII-05}{\it{Таргонський А.Л.}\/} Екстремальні задачі теорії однолистих функцій:
Дис. … канд. физ.-мат. наук. --- Київ, 2005. -- 143 с.

\item\label{kovalev}{\it{Ковалев Л.В.}\/} К задаче об экстремальном разбиении со свободными полюсами на окружности. // Дальневосточный матем. сборник. -- 1996. -- 2. -- С.~96~--~98.

\item\label{Bakhtin}{\it{Бахтин А.К., Денега И.В.}\/} Некоторые оценки функционалов для $N$-лучевых систем точек // Теорія наближення функцій та суміжні питання / Зб. праць Ін-ту матем. НАН України. -- К.: Ін-т матем. НАН України, 2011. -- Т.8, №1. -- C. 12 -- 21.

\end{enumerate}

\end{document}